\documentclass[reqno,10pt]{amsart}

\usepackage{amsmath,amssymb}
\usepackage[dvips]{graphics}
\usepackage{epsfig}

\newtheorem{theorem}{Theorem}[section]

\newtheorem{conjecture}[theorem]{Conjecture}

\newtheorem{definition}[theorem]{Definition}

\newtheorem{lemma}[theorem]{Lemma}

\numberwithin{equation}{section}

\theoremstyle{definition}

\newcommand\ben{\begin{enumerate}}
\newcommand\een{\end{enumerate}}

\newcommand\be{\begin{equation}}
\newcommand\ee{\end{equation}}
\newcommand\benn{\begin{equation*}}
\newcommand\eenn{\end{equation*}}
\newcommand\bea{\begin{eqnarray}}
\newcommand\eea{\end{eqnarray}}
\newcommand\beann{\begin{eqnarray*}}
\newcommand\eeann{\end{eqnarray*}}

\renewcommand{\mod}{\;\operatorname{mod}}

\numberwithin{equation}{section}

\begin{document}

\title{Structure Theorem for $\left( d,g,h\right) $-Maps}

\author{A. V. Kontorovich}
\address{Mathematics Department of Princeton University, Princeton, NJ, $08544$ USA}
\curraddr{Mathematics Department of Columbia University, New York,
NY, $10027$ USA}
\email{alexk@alumni.princeton.edu}
\author{Ya. G. Sinai}
\address{Mathematics Department of Princeton University, Princeton, NJ, $08544$ USA}
\email{sinai@math.princeton.edu}

\begin{abstract}
The $\left( 3x+1\right) $-Map, $T$, acts on the set, $\Pi $, of
positive
integers not divisible by 2 or 3. It is defined by $T\left( x\right) =\frac{%
3x+1}{2^{k}}$, where $k$ is the largest integer for which $T\left(
x\right) $ is an integer. The $\left( 3x+1\right) $-Conjecture asks
if for every $x\in \Pi $ there exists an integer, $n$, such that
$T^{n}\left( x\right) =1$. The Statistical $\left( 3x+1\right)
$-Conjecture asks the same question, except for a subset of $\Pi $
of density 1. The Structure Theorem proven in \cite{sinai} shows
that infinity is in a sense a repelling point, giving some reasons
to expect that the $\left( 3x+1\right) $-Conjecture may be true. In
this paper, we present the analogous theorem for some
generalizations of the $\left( 3x+1\right) $-Map, and expand on the
consequences derived in \cite{sinai}. The generalizations we
consider are determined by positive coprime integers, $d$ and $g$,
with $g>d\geq 2$, and a periodic function, $h\left(
x\right) $. The map $T$ is defined by the formula $T\left( x\right) =\frac{%
gx+h\left( gx\right) }{d^{k}}$, where $k$ is again the largest
integer for which $T\left( x\right) $ is an integer. We prove an
analogous Structure Theorem for $\left( d,g,h\right) $-Maps, and
that the probability distribution corresponding to the density
converges to the Wiener measure with the drift $\log
g-\frac{d}{d-1}\log d$ and positive diffusion constant. This shows
that it is natural to expect that typical trajectores return to the
origin if $\log g-\frac{d}{d-1}\log d<0$ and escape to infinity
otherwise.
\end{abstract}

\keywords{3x+1 Problem, 3n+1 Problem, Collatz Conjecture, Structure
Theorem, ($d$,$g$,$h$)-Maps, Brownian Motion.}

\maketitle

\section{Introduction}

\subsection{\strut The $\left( 3x+1\right) $-Map and $\left( 3x+1\right) $%
-Conjecture}

Recall the definition of the $\left( 3x+1\right) $-Map, (see \cite{lagar}).
Take an integer $x>0$, with $x\ $odd. Then $3x+1$ divides $2$, so we can
find a unique $k>0$ such that $y=\frac{3x+1}{2^{k}}$ is again odd. In this
way, we get a mapping $T:x\longmapsto y$ defined on the set $\Pi $ of
strictly positive numbers not divisible by $2$ or $3$. Write $\Pi =6\Bbb{Z}%
^{+}+E$, where $E=\left\{ 1,5\right\} ,$ is the set of possible congruence
classes modulo 6.

For every integer, $x$, with $0<x<2^{60}$, a computer has checked
that enough iterations of the $\left( 3x+1\right) $-Map eventually
send $x$ to 1 (see \cite{lagar}). The natural conjecture asks if the
same statement holds for all $x\in \Pi $:

\begin{conjecture}[$\left( 3x+1\right) $-Conjecture]
For every $x\in \Pi $, there is an integer $n$, such that $T^{n}\left(
x\right) =1$.
\end{conjecture}

The Statistical $\left( 3x+1\right) $-Conjecture asks the same question,
except for a subset of $\Pi $ of density 1.

\strut

For every $x$, we can associate a value, which is the $k$ used in the
definition of $T$. When we apply $T$ repeatedly, we get a set of $k$ values,
called the path of $x$. We shall call the ordered set of positive integers, $%
\left( k_{1},...,k_{m}\right) $, the ``$m$-path of $x$,'' denoted by $\gamma
_{m}\left( x\right) $, if these are the $k$ values that appear in $m$
repeated iterations of $T$.

E.g. $T\left( 17\right) =\frac{3\cdot 17+1}{2^{2}}=13$, so $k=2,$ and $%
\gamma _{1}\left( 17\right) =\left( 2\right) $. $T^{2}\left( 17\right)
=T\left( 13\right) =\frac{3\cdot 13+1}{2^{3}}=5$, so here $k=3$, and thus $%
\gamma _{1}\left( 13\right) =\left( 3\right) $, and $\gamma _{2}\left(
17\right) =\left( 2,3\right) $.

Assume that we are given an $m$-path, $\left( k_{1},...,k_{m}\right) $. We
can ask the following question: what is the set of $x\in \Pi $ for which $%
\gamma _{m}\left( x\right) =\left( k_{1},...,k_{m}\right) $?

The answer is given by the so-called Structure Theorem, proven in
\cite {sinai}. The theorem states that if $x\in \Pi $ has $\gamma
_{m}\left( x\right) =\left( k_{1},...,k_{m}\right) ,$ then the next
value in $\Pi $ which will have the same $m$-path and congruence
class modulo 6 is $x+6\cdot
2^{k_{1}+...+k_{m}}$. In other words, there is some first $x\in \Pi =6\Bbb{Z}%
^{+}+E$, call it $x_{0}$, which has $\gamma _{m}\left( x\right) =\left(
k_{1},...,k_{m}\right) $. Writing $x_{0}=6\cdot q+\varepsilon $, with $%
\varepsilon \in E$, we get all $x$ with the same $\varepsilon $ and $m$-path
from the sequence $x_{p}=6\left( 2^{k_{1}+...+k_{m}}p+q\right) +\varepsilon $%
. The theorem tells us how to solve uniquely for $q$ given the $m$-path and $%
\varepsilon $, and shows that $q<2^{k_{1}+...+k_{m}}$, so the representation
of $x_{p}$ is unique.

E.g. Let $k_{1}=2$, $k_{2}=3$, and $\varepsilon =5$. Then $x_{0}=17=6\left(
2^{5}\cdot 0+2\right) +5$, and we know that $\gamma _{2}\left( 17\right)
=\left( 2,3\right) $. Look at $x_{1}=6\left( 2^{5}\cdot 1+2\right) +5=209$: $%
T\left( 209\right) =157$, with $k_{1}=2$, and $T^{2}\left( 209\right) =59$
with $k_{2}=3,$ so $\gamma _{2}\left( 209\right) =\left( 2,3\right) $. We
can verify that there are no elements of $\Pi $ between 18 and 208 that are
congruent 5 modulo 6 and have the $2$-path $\left( 2,3\right) $.

Moreover, the Structure Theorem tells us that if the image of $x_{0}$ is $%
y_{0}=T^{m}\left( x_{0}\right) =6\cdot r+\delta $, with $\delta \in E$
(since $y_{0}$ is also in $\Pi $), then we get the next image by adding $%
6\cdot 3^{m}$. In other words, if $y_{p}$ is the image of $x_{p}$, then $%
y_{p}=T^{m}\left( x_{p}\right) =6\left( 3^{m}p+r\right) +\delta $. The
theorem also solves explicitly for $r$ and $\delta $ given the $m$-path and $%
\varepsilon $, and finds that $r<3^{m}$.

\strut E.g. $T^{2}\left( 17\right) =5=6\left( 3^{2}\cdot 0+0\right) +5$, and
$T^{2}\left( 209\right) =59=6\left( 3^{2}\cdot 1+0\right) +5$.

The Structure Theorem also shows that infinity is in a sense a repelling
point. This gives some reasons to expect that the $\left( 3x+1\right) $%
-Conjecture may be true.

In this paper, we present the analogous theorem for some
generalizations of the $\left( 3x+1\right) $-Map, and expand on the
consequences derived in \cite{sinai}.

\subsection{\strut The $\left( d,g,h\right) $-Maps and $\left( d,g,h\right) $%
-Problem}

\strut The generalizations we consider are a particular case of maps
proposed in \cite{RandF}. They are determined by positive coprime
integers, $d$ and $g$, with $g>d\geq 2$, and a periodic function,
$h\left( x\right) $, satisfying:

\begin{enumerate}
\item  $h\left( x+d\right) =h\left( x\right) $,

\item  $x+h\left( x\right) \equiv 0\left( \mod d\right) $,

\item  $0<\left| h\left( x\right) \right| <g$ for all $x$ not divisible by $%
d $.
\end{enumerate}

The map $T$ is defined by the formula
\[
T\left( x\right) =\frac{gx+h\left( gx\right) }{d^{k}},
\]
where $k$ is uniquely chosen so that the result is not divisible by $d$.
Property 2 of $h$ guarantees $k\geq 1$. The natural domain of this map is
the set $\Pi $ of positive integers not divisible by $d$ and $g$. Let $E$ be
the set of integers between $1$ and $dg$ that divide neither $d$ nor $g$, so
we can write $\Pi =dg\Bbb{Z}^{+}+E$. The size of $E$ can easily be
calculated: $\left| E\right| =\left( d-1\right) \left( g-1\right) $.

In the same way as before, we have $m$-paths, which are the values of $k$
that appear in iterations of $T$, and we again denote them by $\gamma
_{m}\left( x\right) $.

The original problem corresponds to $g=3$, $d=2$, and $h\left( 1\right) =1$.
The $\left( 3x-1\right) $-problem corresponds to $g=3$, $d=2$, and $h\left(
1\right) =-1$. The $\left( 5x+1\right) $-problem corresponds to $g=5$, $d=2$%
, and $h\left( 1\right) =1$, and so on.

\strut

The Structure Theorem for $\left( d,g,h\right) $-Maps will be slightly
different, in that given an $m$-path, $\left( k_{1},...,k_{m}\right) $, and
congruence class, $\varepsilon $, modulo $dg$, we do not have a unique $%
x_{0} $. Instead, we have $\left( d-1\right) ^{m}$ values of what was $x_{0}$
in the original case, which we will denote by $x_{0}^{\left( i\right) }$,
with $i=1,...,\left( d-1\right) ^{m}$. Each of these can be written as $%
x_{0}^{\left( i\right) }=dg\cdot q^{\left( i\right) }+\varepsilon $, with $%
q^{\left( i\right) }<d^{k_{1}+...+k_{m}}$. Then we get every $x$ with the
given $m$-path by adding $dg\cdot d^{k_{1}+...+k_{m}}$. In other words,
letting
\[
x_{p}^{\left( i\right) }=dg\left( d^{k_{1}+...+k_{m}}p+q^{\left( i\right)
}\right) +\varepsilon ,
\]
we get every $x\in \Pi $ with $\gamma _{m}\left( x\right) =\left(
k_{1},...,k_{m}\right) $ and $x\equiv \varepsilon \left( \mod
dg\right) $ in the set $\left\{ x_{p}^{\left( i\right) }\right\}
_{p\geq 0,1\leq i\leq \left( d-1\right) ^{m}}$.

\strut Here is the precise formulation of the Structure Theorem for $\left(
d,g,h\right) $-Maps.

\begin{theorem}[Structure Theorem]
Given an $m$-path, $\left( k_{1},...,k_{m}\right) $, and $\varepsilon \in E$%
, let $k=k_{1}+...+k_{m}$. Then there exist $\left( d-1\right) ^{m}$
triples, $\left( q^{\left( i\right) },r^{\left( i\right) },\delta ^{\left(
i\right) }\right) $, $i=1,...,\left( d-1\right) ^{m}$, with $0\leq q^{\left(
i\right) }<d^{k}$, $0\leq r^{\left( i\right) }<g^{m}$, and $\delta ^{\left(
i\right) }\in E$, such that
\[
\left\{ x\in \Pi :x\equiv \varepsilon \left( \mod dg\right)
,\;\gamma _{m}\left( x\right) =\left( k_{1},...,k_{m}\right)
\right\} =\left\{ dg\left( d^{k}p+q^{\left( i\right) }\right)
+\varepsilon \right\} _{p\geq 0,1\leq i\leq \left( d-1\right) ^{m}}.
\]
Moreover, $T^{m}\left( dg\left( d^{k}p+q^{\left( i\right) }\right)
+\varepsilon \right) =dg\left( g^{m}p+r^{\left( i\right) }\right) +\delta
^{\left( i\right) }$.
\end{theorem}

The proof of the theorem is given in the next section.

In section 3, we prove that the probability distribution corresponding to
the density converges to the Wiener measure with the drift $\log g-\frac{d}{%
d-1}\log d$ and positive diffusion constant. This shows that it is natural
to expect that typical trajectories return to the origin if $\log g-\frac{d}{%
d-1}\log d<0$ and escape to infinity otherwise. This question is discussed
in more detail in section 4.

\section{Proof of the Structure Theorem}

The proof goes by induction on $m$. At each stage, we assume $x$ has the
given $m$-path and modulo class, and write $x=dg\left( d^{k}p+q\right)
+\varepsilon $ and $y=T^{m}\left( x\right) =dg\left( g^{m}s+r\right) +\delta
$. This can be done for any number, since we are simply writing out the
modulo classes. After some algebra, we come to some equation for the
triplets $\left( q,r,\delta \right) $, and show that it has $\left(
d-1\right) ^{m}$ solutions.

\subsection{\strut Case $m=1$}

Say we are given a $1$-path, $\left( k\right) $, and let us take an $%
\varepsilon \in E$. Write $x=dg\cdot t+\varepsilon $, and assume that $x$
has the $1$-path, $\left( k\right) $. One can further break $t$ into the
form: $t=d^{k}p+q$, with $0\leq q<d^{k}$. Let $y=T\left( x\right) $, so by
our assumption, $d^{k}y=gx+h\left( gx\right) $. By periodicity, $h\left(
gx\right) =h\left( g\varepsilon \right) $, so since $\varepsilon $ is fixed,
$h$ does not depend on $x$, and is fixed. Thus we will write just $h$ for $%
h\left( gx\right) $ from now on. Since $y\in \Pi $, we can write $y=dg\cdot
t^{\prime }+\delta $ for some $\delta \in E$, and expand $t^{\prime }=g\cdot
s+r$, for $0\leq r<g$. The first step of our analysis is to show that $s=p$.
We write $gx+h=d^{k}y$, and substitute for $x,\;y,\;t,\;$and $t^{\prime }$:
\[
g\left( dg\cdot \left( d^{k}p+q\right) +\varepsilon \right) +h=d^{k}\left(
dg\cdot \left( g\cdot s+r\right) +\delta \right) .
\]
We expand this to see:
\begin{equation}
g^{2}d^{k+1}\cdot p+\left( dg^{2}q+g\varepsilon +h\right) =g^{2}d^{k+1}\cdot
s+\left( d^{k+1}gr+d^{k}\delta \right) .  \label{parenth}
\end{equation}
Next, we apply the following simple Lemma.

\begin{lemma}
If $a\cdot b+c=a\cdot b^{\prime }+c^{\prime }$ with $0\leq c,c^{\prime }<a$,
then $b=b^{\prime }$ and $c=c^{\prime }$.%
\endproof%
\end{lemma}

\strut

To apply the lemma (with $a=g^{2}d^{k+1}$), we need to show that the parts
in parentheses on both sides of $\left( \ref{parenth}\right) $ are contained
in $\left[ 0,g^{2}d^{k+1}-1\right] $. We will derive upper and lower bounds
for the left side, and leave similar calculations for the right side to the
reader.

Consider the lower bound of the left side. Since $q\geq 0,\;\varepsilon \geq
1$ and $h\geq -g+1$ (by Condition 3), we have that
\[
dg^{2}\cdot q+g\varepsilon +h\geq g\cdot 1+\left( -g+1\right) =1,
\]
and thus is positive.

For the upper bound of the left side, we notice that $q\leq d^{k}-1$, $%
\varepsilon \leq dg-1$ (since $\varepsilon \in E$) and $h\leq g-1$. So
\begin{eqnarray*}
dg^{2}\cdot q+g\varepsilon +h &\leq &g^{2}d\cdot \left( d^{k}-1\right)
+g\left( dg-1\right) +\left( g-1\right) \\
&=&g^{2}d^{k+1}-1.
\end{eqnarray*}

\strut

\begin{itemize}
\item  The Lemma gives us that $p=s,$ and from now on we write just $p$. We
want to characterize $q,$ $r$ and $\delta $, showing that they are
independent of $p$.
\end{itemize}

\strut

To continue, we recall that the Lemma implies that the parts in parentheses
of $\left( \ref{parenth}\right) $ also concur. So:

\begin{equation}
g^{2}d\cdot q+g\varepsilon +h=d^{k}gd\cdot r+d^{k}\delta .  \label{delt find}
\end{equation}
The next step is to break $\delta $ into $\delta =\delta ^{\prime }g+\delta
^{\prime \prime }$, with $0\leq \delta ^{\prime \prime }<g$. Since $\delta
\in E$, we have $\delta <dg$, implying $0\leq \delta ^{\prime }<d$. We now
look at $\left( \ref{delt find}\right) $ modulo $g$ to solve for $\delta
^{\prime \prime }$:
\begin{equation}
d^{k}\delta ^{\prime \prime }=h\left( \mod g\right) . \label{delt''}
\end{equation}
Since $g$ and $d$ are relatively prime, $d^{k}$ has a multiplicative inverse
in $\left( \Bbb{Z}\setminus g\Bbb{Z}\right) ^{*}$, meaning $\delta ^{\prime
\prime }$ is uniquely determined. Exactly one of the $d$ possible values of $%
\delta ^{\prime }$ will make $\delta =\delta ^{\prime }g+\delta ^{\prime
\prime }$ divisible by $d$, and we throw this value away since $\delta \in E$%
.

\strut

This leaves us with $d-1$ possible values for $\delta $, which we denote by $%
\delta ^{\left( 1\right) },\delta ^{\left( 2\right) },...,\delta ^{\left(
d-1\right) }$. It suffices to solve $\left( \ref{delt find}\right) $
uniquely for $q^{\left( i\right) }$ and $r^{\left( i\right) }\ $given $%
\delta ^{\left( i\right) }$.

\strut

\strut Now we assume we have fixed $\delta ^{\left( i\right) }$, and
rearrange $\left( \ref{delt find}\right) $, adding a superscript to $q$ and $%
r$ to correspond to $\delta $:

\[
g\cdot q^{\left( i\right) }-d^{k}r^{\left( i\right) }=\frac{d^{k}\delta
^{\left( i\right) }-g\varepsilon -h}{dg}=v.
\]
Everything on the right hand side is known, so $v$ is now just an integer
(and independent of $p$). We solve for $q^{\left( i\right) }$ and $r^{\left(
i\right) }$ by applying the Chinese Remainder Theorem to the equation $%
g\cdot a-d^{k}b=1$, then setting $q^{\left( i\right) }=v\cdot
a\left( \mod d^{k}\right) $ and $r^{\left( i\right) }=v\cdot b\left(
\mod g\right) $.

\strut

\begin{itemize}
\item  Having found the triplets $\left( q^{\left( i\right) },r^{\left(
i\right) },\delta ^{\left( i\right) }\right) $, we are done with the case $%
m=1$.
\end{itemize}

\strut

Summarizing the first step of the induction, we pick some $\varepsilon \in E$%
, assume $x\in \Pi $ is of the form $dg\cdot t+\varepsilon $, and write $%
t=d^{k}p+q$. Under the same assumptions for the image, $y=T\left( x\right) $%
, we write $y=dg\cdot t^{\prime }+\delta $ and $t^{\prime }=gp+r$. We find
that $\delta $ is unique modulo $g$, and there are $d-1$ values, $\delta
^{\left( 1\right) },...,\delta ^{\left( d-1\right) }$, which $\delta \in E$
may take. For each one, we solve for $q^{\left( i\right) }$ and $r^{\left(
i\right) }$. All of the calculations depend only on $k$ and $\varepsilon $.

\subsection{\strut Induction on $m>1$}

For $m>1$, the induction goes as follows. To know which $x$ have a given $m$%
-path, $\left( k_{1},...,k_{m}\right) $, we first assume we know the answer
for the $\left( m-1\right) $-path, $\left( k_{1},...,k_{m-1}\right) $.

Let $k=k_{1}+k_{2}+...+k_{m-1}$, and assume by the induction hypothesis that
there are $\left( d-1\right) ^{m-1}$ values for the triplet $\left(
q_{m-1},r_{m-1},\delta _{m-1}\right) $ which satisfy our equations. Fix one
such triplet, pick any integer, $p_{m-1}$, and set $x=dg\left(
d^{k}p_{m-1}+q_{m-1}\right) +\varepsilon ,$ and $y=dg\left(
g^{m-1}p_{m-1}+r_{m-1}\right) +\delta _{m-1}$. Then we have $\gamma
_{m}\left( x\right) =\left( k_{1},...,k_{m-1}\right) $, and $y=T^{m-1}\left(
x\right) $. Here we write $p_{m-1}$ instead of just $p$ to distinguish from
the $p$ we will have in the next paragraph. The triplet $\left(
q_{m-1},r_{m-1},\delta _{m-1}\right) $ is still gotten independently of $%
p_{m-1}$.

We can alternatively break $x$ into $x=dg\left(
d^{k+k_{m}}p_{m}+q_{m}\right) +\varepsilon $ for some $q_{m}<d^{k+k_{m}}$
and also write $z=T^{m}\left( x\right) =T\left( y\right) =dg\cdot t+\delta
_{m}$, with $t=g^{m}s+r_{m}$. The key idea is to find the $d-1$ possible
values for $\delta _{m}\in E$, and with each we solve for the corresponding $%
q_{m}$ and $r_{m}$, knowing $q_{m-1}$, $r_{m-1}$, and $\delta _{m-1}$. We
will again see that $p_{m}=s$ and that $\left( q,r,\delta \right) $ do not
depend on this value.

\strut

Since $z=T\left( y\right) $, by assumption, we have $d^{k_{m}}z=gy+h\left(
gy\right) $, (again let $h=h\left( gy\right) =h\left( g\delta _{m-1}\right) $%
) which expands to:
\begin{equation}
d^{k_{m}+1}g^{m+1}s+d^{k_{m}+1}gr_{m}+d^{k_{m}}\delta
_{m}=dg^{m+1}p_{m-1}+g^{2}dr_{m-1}+g\delta _{m-1}+h.  \label{return}
\end{equation}
Remembering the two expressions for $x,$ and setting $%
p_{m}=d^{k_{m}}p_{1}+p_{2}$ (with $0\leq p_{2}<d^{k_{m}}$), we write:
\begin{eqnarray*}
d^{k+k_{m}}p_{m}+q_{m} &=&\frac{x-\varepsilon }{dg}=d^{k}p_{m-1}+q_{m-1} \\
&=&d^{k+k_{m}}p_{1}+d^{k}p_{2}+q_{m-1}.
\end{eqnarray*}
We easily see that $0\leq d^{k}p_{2}+q_{m-1}<d^{k+k_{m}}$, so we again use
the Lemma to find:
\begin{eqnarray}
p_{m} &=&p_{1},  \label{p_m=p_1} \\
q_{m} &=&d^{k}p_{2}+q_{m-1}.  \label{q_m=}
\end{eqnarray}

\strut Returning to $\left( \ref{return}\right) $, we expand:
\[
d^{k_{m}+1}g^{m+1}s+\left( d^{k_{m}+1}gr_{m}+d^{k_{m}}\delta _{m}\right)
=d^{k_{m}+1}g^{m+1}p_{1}+\left( dg^{m+1}p_{2}+g^{2}dr_{m-1}+g\delta
_{m-1}+h\right) .
\]

\strut Following the same techniques as before, we bound the parts in
parentheses on both sides between zero and $d^{k_{m}+1}g^{m+1}$, and apply
the Lemma. This gives us that $p_{m}=p_{1}=s$, and that
\begin{equation}
d^{k_{m}+1}gr_{m}+d^{k_{m}}\delta _{m}=dg^{m+1}p_{2}+g^{2}dr_{m-1}+g\delta
_{m-1}+h.  \label{close}
\end{equation}
Again looking modulo $g$ and setting $\delta _{m}=\delta ^{\prime }g+\delta
^{\prime \prime }$, we solve:
\[
\delta ^{\prime \prime }\equiv g\delta _{m-1}+h\left( \mod g\right)
,
\]
which again gives us $d$ choices for $\delta ^{\prime }$, one of which we
throw out because $\delta _{m}\in E$. Rearranging $\left( \ref{close}\right)
$, we get:
\[
g^{m}p_{2}-d^{k_{m}}r_{m}=\frac{d^{k_{m}}\delta _{m}-dg^{2}r_{m-1}-g\delta
_{m-1}-h}{dg}=v
\]

\strut From here, we solve $g^{m}a-d^{k_{m}}b=1$ and set
$p_{2}=a\cdot v\left( \mod d^{k_{m}}\right) $ and $r_{m}=b\cdot
v\left( \mod g^{m}\right) $, so $q_{m}=d^{k}p_{2}+q_{m-1}$. We have
$\left( d-1\right) $ values of $\left( q_{m},r_{m},\delta
_{m}\right) $ derived from $\left( d-1\right) ^{m-1}$ values of
$\left( q_{m-1},r_{m-1},\delta _{m-1}\right) $, so there are a total
of $\left( d-1\right) ^{m}$ triplets, consistent with the induction
hypothesis. Now everything in the triplet $\left(
q_{m},r_{m},\delta _{m}\right) $ is defined, and we are done.%
\endproof%

\section{Brownian Motion of $\left( d,g,h\right) $-Paths}

\strut In \cite{FMMT},\cite{LandW}, it is assumed that the $\left(
3x+1\right) $-Map behaves as a geometric Brownian motion, and a
stochastic model is built from which other conjectures relating to
the problem are derived. Here, we prove that the generalized $\left(
d,g,h\right) $-Maps do indeed have this behavior.

In order to consider sample $\left( d,g,h\right) $-paths, we must first
establish a version of a probability measure on $\Bbb{Z}^{+}$. The only
natural way to do this is through density:

\begin{definition}
For $A\subset \Bbb{Z}^{+}$, define
\begin{equation}
P\left( A\right) =\lim_{n\rightarrow \infty }\frac{\left| A\cap \left[
1,n\right] \cap \Pi \right| }{\left| \left[ 1,n\right] \cap \Pi \right| }%
=\lim_{n\rightarrow \infty }\frac{\left| A\cap \left[ 1,n\right] \cap \Pi
\right| }{n}\cdot \frac{dg}{\left| E\right| },  \label{WorkPlease}
\end{equation}
provided the limit exists.
\end{definition}

A nice consequence of the Structure Theorem is that if we want to consider
the set of $x$ that follow a certain $m$-path, they all fall in one of
several arithmetic progressions, and so these sets have a density.

\strut

Partition the interval $\left[ 0,1\right] $ by:\ $0=t_{0}<t_{1}<...<t_{r}=1$%
. Fix $m$ and let $m_{i}=\left\lfloor t_{i}m\right\rfloor $. For any $x$,
let $x_{i}=T^{m_{i}}\left( x\right) $.

\begin{theorem}
The properly normalized path $\ln x_{i}$ converges as $m\rightarrow \infty $
to a Brownian Path with drift $\ln g-\frac{d}{d-1}\ln d$. More precisely,
\begin{eqnarray*}
&&\lim_{m\rightarrow \infty }P\left\{ x:a_{i}<\frac{\ln x_{i+1}-\ln
x_{i}-\left( m_{i+1}-m_{i}\right) \left( \ln g-\frac{d}{d-1}\ln d\right) }{%
\sqrt{\frac{d}{\left( d-1\right) ^{2}}m}\ln d}<b_{i},\text{ with }%
i=0,...,r-1\right\} \\
&& \\
&&\;\;\;\;\;\;\;\;\;\;\;\;\;\;\;\;\;\;\;\;\;\;\;\;\;\;\;\;\;\;\;\;\;\;\;\;\;%
\;\;\;\;\;\;\;\;\;\;\;\;\;\;\;\shortparallel \\
&& \\
&&\;\;\;\;\;\;\;\;\;\;\;\;\;\;\;\;\;\;\;\;\;\;\;\;\;\;\;\;\;\;\;\;\;\;\;\;\;%
\int_{a_{0}}^{b_{0}}\int_{a1}^{b_{1}}\cdots \int_{a_{r-1}}^{b_{r-1}}\frac{%
e^{\left( -\frac{1}{2}\sum_{i=0}^{r-1}u_{i}^{2}\right) }}{\left( 2\pi
\right) ^{\frac{r}{2}}}du_{0}du_{1}\cdots du_{r-1}.
\end{eqnarray*}
\end{theorem}

\strut

\proof%
By an extension of the Structure Theorem, we know that $x_{i}=T^{m_{i}}%
\left( x\right) $ can be expressed as $x_{i}=dg\left(
g^{m_{i}}d^{k_{m_{i}+1}+...+k_{m}}p+q_{i}\right) +\delta _{i}$. Then
\begin{equation}
\ln x_{i}=m_{i}\ln g+\left( k_{m_{i}+1}+...+k_{m}\right) \ln d+\ln p+O\left(
1\right) ,  \label{logX}
\end{equation}
and since we are interested in questions about density, $x_{i}$ is large, so
$p$ is large, and thus $O\left( 1\right) $ is non-essential. Then we can
rearrange $\left( \ref{logX}\right) $ to:
\begin{eqnarray*}
\ln x_{i}-m_{i}\ln g-\left( k_{m_{i}+1}+...+k_{m}\right) \ln d &=&\ln p \\
&=&\ln x_{i+1}-m_{i+1}\ln g-\left( k_{m_{i+1}+1}+...+k_{m}\right) \ln d,
\end{eqnarray*}
from which we get:
\begin{equation}
\left( m_{i+1}-m_{i}\right) \frac{d}{d-1}\ln d-\left(
k_{m_{i}+1}+...+k_{m_{i+1}}\right) \ln d=\ln x_{i+1}-\ln x_{i}-\left(
m_{i+1}-m_{i}\right) \left( \ln g-\frac{d}{d-1}\ln d\right) .
\label{diffRelation}
\end{equation}
Since the set of $x_{i}$ consists of precisely $\left( d-1\right) ^{i}$
arithmetic progressions, each with step $dg\cdot d^{k}$ (where $%
k=k_{1}+...+k_{m}$), we use $\left( \ref{WorkPlease}\right) $ to find that
\[
P\left\{ \gamma _{m}\left( x\right) =\left( k_{1},...,k_{m}\right) ,x\equiv
\varepsilon \left( \mod dg\right) \right\} =\frac{1}{dg\cdot d^{k}}%
\frac{dg}{\left| E\right| }\left( d-1\right) ^{m}.
\]
This holds for each $\varepsilon \in E$, so we see that
\begin{eqnarray}
P\left\{ \gamma _{m}\left( x\right) =\left( k_{1},...,k_{m}\right)
\right\} &=&\left| E\right| \cdot P\left\{ \gamma _{m}\left(
x\right) =\left( k_{1},...,k_{m}\right) ,x\equiv \varepsilon \left(
\mod dg\right)
\right\}  \nonumber \\
&=&\frac{\left( d-1\right) ^{m}}{d^{k}}=\prod_{j=1}^{m}\frac{\left(
d-1\right) }{d^{k_{j}}}.  \label{kdist}
\end{eqnarray}
This shows that we can consider the $k_{j}$ as independent identically
distributed random variables, with exponential distribution having the
parameter $\frac{1}{d}$. Thus the expected value,
\begin{eqnarray*}
E\left[ k_{1}+...+k_{m}\right] &=&\sum_{n\geq m}n\cdot P\left\{
k_{1}+...+k_{m}=n\right\} \\
&=&\sum_{n\geq m}n\cdot \sum_{s_{1}+...+s_{m}=n-m\text{, }s_{i}\geq
0}P\left\{ \left( s_{1}+1,...,s_{m}+1\right) \right\} \\
&=&\left( d-1\right) ^{m}\sum_{n\geq m}n\sum_{s_{1}+...+s_{m}=n-m\text{, }%
s_{i}\geq 0}\frac{1}{d^{n}} \\
&=&\left( d-1\right) ^{m}\sum_{n\geq m}\frac{n}{d^{n}}\left(
\begin{array}{c}
n-1 \\
m-1
\end{array}
\right) \\
&=&\frac{d}{d-1}m.
\end{eqnarray*}
Similarly, we can calculate that $Var\left[ k_{1}+...+k_{m}\right] =\frac{d}{%
\left( d-1\right) ^{2}}m$. So by the Central Limit Theorem,
\[
\lim_{m\rightarrow \infty }P\left\{ \frac{k_{1}+...+k_{m}-\frac{d}{d-1}m}{%
\sqrt{\frac{d}{\left( d-1\right) ^{2}}m}}\in \left( a,b\right) \right\}
=\int_{a}^{b}\frac{e^{-\frac{u^{2}}{2}}}{\sqrt{2\pi }}du.
\]
And by $\left( \ref{diffRelation}\right) $, we have that
\begin{eqnarray*}
&&P\left\{ \frac{\ln x_{i+1}-\ln x_{i}-\left( m_{i+1}-m_{i}\right) \left(
\ln g-\frac{d}{d-1}\ln d\right) }{\sqrt{m\cdot \frac{d}{\left( d-1\right)
^{2}}}\ln d}\in \left( a_{i},b_{i}\right) \right\} \\
&& \\
&&\;\;\text{\ \ }\;\;\;\;\;\;\;\;\;\;\;\;\;\;\;\;\;\;\;\;\;\;\;\;\;\;\;\;\;%
\;\;\;\;\shortparallel \\
&& \\
&&\;\;\;P\left\{ \frac{\frac{d}{d-1}\left( m_{i+1}-m_{i}\right) -\left(
k_{m_{i}+1}+...+k_{m_{i+1}}\right) }{\sqrt{\frac{d}{\left( d-1\right) ^{2}}m}%
}\in \left( a_{i},b_{i}\right) \right\} ,
\end{eqnarray*}
which converges exactly as claimed. Since the $k_{i}$ are independent, the
increments, $\ln x_{i+1}-\ln x_{i}$ are as well, and we have the statement
about the convergence of our distributions to the Wiener measure.%
\endproof%

\strut

\section{Asymptotic Behavior of Typical Trajectories}

The previous section proves that the probability distribution corresponding
to the density converges to the Wiener measure with drift $\log g-\frac{d}{%
d-1}\log d$. Since $d$ and $g$ are relatively prime, there are no values of $%
d$ and $g$ for which $\log g-\frac{d}{d-1}\log d=0$, and thus every $\left(
d,g,h\right) $-Map has a non-trivial drift. Therefore, the asymptotic
behavior of typical trajectories depends entirely on the sign of the drift.
When the drift is negative, infinity is a repelling point. In the opposite
case, typical trajectories escape to infinity. For the original $\left(
3x+1\right) $-Map, the drift is $\log 3-2\log 2<0$, and so as a special
case, we get the result found in \cite{sinai}.

In the literature, the stopping time of an integer $x$ is defined as
the first positive integer, $n,$ such that $T^{n}\left( x\right)
<x$. If $n$ does not exist, we say that $x$ has an infinite stopping
time. In \cite {Everett} and \cite{Tarras},\cite{Terras}, it is
independently proven that for the $\left( 3x+1\right) $-Map, the
density of integers with a finite stopping time is 1. This paper
provides another proof of this statement.

\strut

\strut

\textbf{Acknowledgments:}\ The first author thanks L. Kontorovich and S.
Payne for discussions and criticism. The second author thanks the NSF\ for
financial support, grant DMR-9813268.

\strut


\begin{thebibliography}{FMMT}
\bibitem[E]{Everett}  C. J. Everett, \textit{Iteration of the number
theoretic function }$f\left( 2n\right) =n$\textit{, }$f\left( 2n+1\right)
=3n+2$, Adv. Math., \textbf{25 }(1977), 42-45.

\bibitem[FMMT]{FMMT}  Feix, M.R.; Muriel, A.; Merlini, D.; Tartini, R.
\textit{The }$\left( 3x+1\right) /2$ \textit{problem: A Statistical Approach}%
, in: Stochastic Processes, Physics and Geometry II, Locarno 1991. World
Scientific (1995), 289-300.

\bibitem[FR]{RandF}  Feix, M.R. and Rouet, J.L. \textit{The }$\left(
3x+1\right) /2$ \textit{problem and its generalization: a stochastic approach%
},\ Proceedings of the International Conference on the Collatz Problem and
Related Topics (2001).

\bibitem[L]{lagar}  Lagarias, J.C., \textit{The }$3x+1$\textit{\ Problem and
Its Generalizations,} American Mathematical Monthly, Vol. \textbf{92}, Issue
1 (Jan., 1985), 2-23.

\bibitem[LW]{LandW}  Lagarias, J.C., and Weiss, A. \textit{The }$3x+1$%
\textit{\ Problem and:\ Two stochastic models, }Ann. Appl. Prob. $\mathbf{2}$
(1992) 229-261.

\bibitem[S]{sinai}  Sinai, Ya. G., \textit{Statistical }$\left( 3x+1\right) $%
\textit{-Problem,} (2002), preprint.

\bibitem[T76]{Tarras}  R. Terras, \textit{A\ stopping time problem on the
positive integers,} Acta Arith. \textbf{30} (1976), 241-252.

\bibitem[T79]{Terras}  R. Terras, \textit{On the existence of a density,}
Acta Arith. \textbf{35} (1979), 101-102.
\end{thebibliography}
\end{document}